\documentclass[conference]{IEEEtran}
\IEEEoverridecommandlockouts
\usepackage{comment}
\usepackage{cite}
\usepackage{amsmath,amssymb,amsfonts}
\usepackage{algorithmic}
\usepackage{algorithm}
\usepackage{graphicx}
\usepackage{textcomp}
\usepackage{booktabs}
\usepackage{multirow}
\usepackage{xcolor}
\usepackage{subcaption}
\usepackage{hyperref}
\usepackage{xurl}
\def\BibTeX{{\rm B\kern-.05em{\sc i\kern-.025em b}\kern-.08em
    T\kern-.1667em\lower.7ex\hbox{E}\kern-.125emX}}

\newcommand{\dt}{\Delta t}
\newcommand{\feps}{f_{\epsilon}}

\begin{document}

\title{Performance Evaluation of Stabilized Corrections\\ for Mixed Precision Runge--Kutta Methods}

\author{
\IEEEauthorblockN{C\'esar Herrera}
\IEEEauthorblockA{\textit{Department of Mathematics} \\
\textit{Purdue University}\\
West Lafayette, IN, USA \\
herre125@purdue.edu}
 \and
 \IEEEauthorblockN{John Driscoll \& Sigal Gottlieb \\
 \& Zachary J. Grant \& Tej Sai Kakumanu}
 \IEEEauthorblockA{\textit{Mathematics Department, UMass Dartmouth} \\
 \textit{North Dartmouth, MA, USA} \\
 jdriscoll9@umassd.edu \& sgottlieb@umassd.edu \\
 zgrant@umassd.edu \& tkakumanu@umassd.edu \\
 \vspace*{-.6in}}
 \and
 \IEEEauthorblockN{Andrew Christlieb}
 \IEEEauthorblockA{\textit{Department of CMS\&E}\\
 \textit{Michigan State University}\\
 \textit{East Lansing, MI, USA} \\
 christli@msu.edu} 
}
\maketitle

\begin{abstract}
Mixed precision Runge–Kutta methods reduce the cost of the expensive implicit solves in diagonally implicit Runge-–Kutta (DIRK) schemes by evaluating them in low precision, while retaining the accuracy of the scheme for larger time steps. The accuracy lost to the low precision perturbation can be recovered through inexpensive explicit corrections; however, these corrections have an adverse impact on stability. Recently proposed stabilized corrections remedy this by applying a stabilization matrix to the correction step, but their runtime cost has not previously been quantified. In this work, we present a numerical study of the runtime performance of these stabilized corrections. Using spectral semi-discretizations of two nonlinear partial differential equations, the inviscid Burgers' equation and the porous medium equation, we compare uncorrected mixed precision DIRK methods against explicitly corrected and stabilized variants across half, single, double, and quadruple precision pairings, for SDIRK methods of orders two through four. We report convergence, runtime, and speedups, and show that the stabilized corrections improve the accuracy of the mixed precision schemes while preserving substantial runtime savings. All
experiments were performed on an \emph{Intel Xeon Platinum 8480+} CPU with
\emph{Julia version 1.11.4}.
\end{abstract}

\section{Introduction}
The use of mixed precision to accelerate numerical algorithms has become increasingly popular in recent years (see, e.g., \cite{gupta2015, higham2019, haidar2018}), driven both by the runtime benefits of low precision arithmetic and by the proliferation of hardware that natively supports multiple floating-point formats. Such algorithms aim to combine the efficiency of low precision with the accuracy of high precision.

In the context of ordinary differential equations and the time-evolution of partial differential equations, Z. J. Grant proposed in \cite{grant2022} a perturbed Runge--Kutta framework that treats the use of a cheaper, lower-accuracy implicit solve as an additive perturbation of  a Runge--Kutta method, enabling the design of mixed precision methods whose low precision errors are suppressed. The performance of these methods was subsequently evaluated in \cite{burnett2021}, where speedups ranging from roughly $2\times$ to $15\times$ were reported, and their stability was analyzed in \cite{burnett2024}.

The accuracy lost to the low precision perturbation can be recovered by appending inexpensive, explicit high precision corrections to the implicit stages \cite{grant2022}. Although these corrections are effective for sufficiently small time steps, it was numerically shown in \cite{burnett2024} that they shrink the region of linear stability and can introduce instabilities for larger time steps. To address this limitation, stabilized corrections were proposed in \cite{driscoll2026}. While the analysis and numerical experiments in \cite{driscoll2026} show that these new correction approaches improve both accuracy and stability, they do not quantify the associated runtime costs. The purpose of this paper is to evaluate the computational runtimes of these stabilized corrections. When applying these stabilized corrections to mixed-precision DIRK methods of second, third, and fourth order, our numerical experiments show speedups ranging from $5\times$ to $55\times$ when combining half precision with double or quadruple precision.


\section{Mixed precision DIRK Methods \& Corrections}\label{sec:methods}
We consider an initial value problem
\begin{equation}\label{eq:ivp}
  y' = f(y), \qquad y(0) = y_0,
\end{equation}
where $f$ is contractive, which typically results from a spatial
semi-discretization of a PDE. We also assume that the derivative of $f$ is
bounded, i.e.\ $\|f'(y)\| \le L$ for some $L > 0$.

\vspace*{-.02in}

\subsection{Mixed precision DIRK Methods}
A diagonally implicit Runge--Kutta (DIRK) method with $s$ stages advances the
solution as
\begin{align*}
  z^{(i)}   &= z_n + \dt \sum_{j=1}^{i} a_{ij}\, f\!\left(z^{(j)}\right), \\
  z_{n+1}   &= z_n + \dt \sum_{i=1}^{s} b_i\, f\!\left(z^{(i)}\right), 
\end{align*}
where $\mathbf{A} = (a_{ij})$ is the matrix of stage coefficients and
$\mathbf{b} = (b_i)$ is the vector of weights.
Following \cite{grant2022,driscoll2026}, we make the implicit stage cheaper by
replacing $f$ with a function $\feps$ in the diagonal (implicit) terms:
\begin{subequations}\label{eq:mp}
\begin{align}
  y^{(i)} &= y_n + \dt\!\left[\sum_{j=1}^{i-1} a_{ij} f\!\left(y^{(j)}\right)
              + a_{ii}\, \feps\!\left(y^{(i)}\right)\right], \label{eq:mp-stage}\\
  y_{n+1} &= y_n + \dt \sum_{i=1}^{s} b_i\, f\!\left(y^{(i)}\right),
            \label{eq:mp-update}
\end{align}
\end{subequations}
where $\epsilon$ represents the precision level of the lower-precision evaluation
$\feps$. The order conditions derived in \cite{grant2022} show that the global
error from replacing $f$ with $\feps$ contributes a term
$\mathcal{O}(\epsilon\, \dt)$, which dominates the scheme error once $\dt$ is
sufficiently small.
We study three DIRK methods along with their mixed precision versions. The first
is the second-order one-stage SDIRK2, i.e.\ the implicit midpoint rule (IMR), with
\begin{equation*}\label{eq:sdirk2-coeff}
  \mathbf{A} = \begin{pmatrix} \tfrac{1}{2} \end{pmatrix}, \qquad
  \mathbf{b} = \begin{pmatrix} 1 \end{pmatrix}.
\end{equation*}
The second is the third-order two-stage SDIRK3 \cite{norsett1974}, with
$\gamma = \tfrac{\sqrt{3}+3}{6}$ and
\begin{equation*}\label{eq:sdirk3-coeff}
  \mathbf{A} = \begin{pmatrix} \gamma & 0 \\[2pt] 1 - 2\gamma & \gamma \end{pmatrix}, \qquad
  \mathbf{b} = \begin{pmatrix} \tfrac{1}{2} \\[2pt] \tfrac{1}{2} \end{pmatrix}.
\end{equation*}
The third is the fourth-order three-stage SDIRK4 \cite{crouzeix1979}, with
$\alpha = \tfrac{2}{\sqrt{3}}\cos\!\left(\tfrac{\pi}{18}\right)$ and
\begin{equation*}\label{eq:sdirk4-coeff}
  \mathbf{A} = \begin{pmatrix}
        \tfrac{1+\alpha}{2} & 0 & 0 \\[4pt]
        -\tfrac{\alpha}{2} & \tfrac{1+\alpha}{2} & 0 \\[4pt]
        1+\alpha & -(1+2\alpha) & \tfrac{1+\alpha}{2}
      \end{pmatrix}, \quad
  \mathbf{b} = \begin{pmatrix}
        \tfrac{1}{6\alpha^{2}} \\[4pt]
        1 - \tfrac{1}{3\alpha^{2}} \\[4pt]
        \tfrac{1}{6\alpha^{2}}
      \end{pmatrix}.
\end{equation*}

\subsection{Explicit Corrections}
To recover accuracy, explicit corrections were introduced in \cite{grant2022}.
For any order-$p$ method \eqref{eq:mp}, we consider the explicit correction
method with $p-1$ correction steps. For each stage $i = 1, \dots, s$, we compute
\begin{align*}
  y^{(i)}_{[0]} &= y_n + \dt\!\left(\sum_{j=1}^{i-1} a_{ij}\, f\!\left(y^{(j)}_{[p-1]}\right)
                   + a_{ii}\, \feps\!\left(y^{(i)}_{[0]}\right)\right), \\
  y^{(i)}_{[k]} &= y_n + \dt\!\left(\sum_{j=1}^{i-1} a_{ij}\, f\!\left(y^{(j)}_{[p-1]}\right)
                   + a_{ii}\, f\!\left(y^{(i)}_{[k-1]}\right)\right), 
\end{align*}
for $k=1, \dots, p-1$, followed by the update
\begin{align*}
  y_{n+1} &= y_n + \dt \sum_{i=1}^{s} b_i\, f\!\left(y^{(i)}_{[p-1]}\right). 
\end{align*}
It was shown in \cite{grant2022} that, asymptotically, each correction
mitigates the perturbation error by a factor of $\dt$. However, as investigated in
\cite{burnett2021,burnett2024}, these corrections shrink the region of linear
stability and may introduce significant instability for larger values of $\dt$.
Following \cite{driscoll2026}, viewing a single implicit stage in the generic form
$y = y_{\mathrm{exp}} + \alpha \dt\, f(y)$, the explicit correction is the
fixed-point iteration
\begin{equation}\label{eq:exp-fp}
  y_{[k+1]} = y_{\mathrm{exp}} + \alpha \dt\, f\!\left(y_{[k]}\right),
\end{equation}
where $\alpha = a_{ii}$ is the diagonal coefficient of the stage. This
iteration converges when $\alpha \dt L < 1$, but may diverge, and the scheme
may become unstable, when $\alpha \dt L \geq  1$.

\subsection{Stabilized Corrections}
To retain the accuracy gain of \eqref{eq:exp-fp} while restoring stability,
\cite{driscoll2026} adds a stabilization term:
\begin{equation}\label{eq:stab-corr}
  y_{[k+1]} = \underbrace{y_{\mathrm{exp}} + \alpha \dt\, f\!\left(y_{[k]}\right)}_{\text{explicit correction}}
            + \underbrace{\alpha \dt\, J\!\left(y_{[k+1]} - y_{[k]}\right)}_{\text{stabilization}}.
\end{equation}
Letting $r_{[k]} = y_{\mathrm{exp}} + \alpha \dt\, f(y_{[k]}) - y_{[k]}$ be the
residual of the implicit stage, \eqref{eq:stab-corr} is equivalent to the update
\begin{equation*}\label{eq:stab-update}
  y_{[k+1]} = y_{[k]} + \Phi\, r_{[k]}, \qquad
  \Phi = \left(I - \alpha \dt\, J\right)^{-1}.
\end{equation*}
 Setting $\Phi = I$ (i.e.\ $J = 0$)
recovers the explicit correction. As presented in \cite{driscoll2026}, there are
several possible choices of $\Phi$. Here we focus on evaluating the performance of
static approaches, in which the matrix $\Phi$ is computed once at the beginning of
the simulation and reused at every time step and iteration. We consider the
following two choices.

\noindent {\bf a) Jacobian-based approach $\Phi_J$:} We freeze the Jacobian $J_0 = f'(y_0)$
at the initial value $y_0$, and compute $ \Phi = \Phi_J = \left(I - \alpha \dt\, J_0\right)^{-1}$ 
once in high precision. $\Phi_J$ is then  re-used throughout the simulation.

\noindent{\bf b) Differential-operator approach $\Phi_\mathrm{EIN}$:} We replace the Jacobian by the
dominant differential operator $\mathcal{L}$ of $f$, {\em e.g.}
$\mathcal{L} = -D_x$ for an advective flux $f(u) = -D_x F(u)$, and
$\mathcal{L} = D_{xx}$ for a diffusion-type operator $f(u) = D_{xx}(u^m)$.
We then compute 
$\Phi = \Phi_{\mathrm{EIN}} = \left(I - \alpha \dt\, \mathcal{L}\right)^{-1}$
in high precision, and re-use it throughout the simulation.
This approach is  inspired by the explicit--implicit--null (EIN) splitting.

\section{Implementation Details}\label{sec:impl}
The methods of Section~\ref{sec:methods} are implemented in the Julia language
using the \texttt{Float16}, \texttt{Float32}, \texttt{Float64}, and
\texttt{Float128} types (software quadruple precision via \texttt{Quadmath.jl}).
The implicit stages are solved with Newton iterations using a stopping tolerance
of $10\,\epsilon_{\text{mach}}$, where $\epsilon_{\text{mach}}$ is the machine
epsilon of the working precision. The maximum number of iterations is set to
$20$, regardless of precision.

The spatial operators are Fourier spectral differentiation matrices, constructed separately in each working precision. For the mixed-precision methods, the stabilization matrices $\Phi_J$ and $\Phi_{\mathrm{EIN}}$ are formed and factored once, in high precision, at the initial time step.

A reference solution for each test problem is computed with the fourth-order
explicit method RK4 in quadruple precision, using a small time step. Runtimes are
reported as the median of repeated timed runs of the time-stepping routine.
We denote a precision pairing as \texttt{full/reduced}, e.g.\ \texttt{64/32} for a
double/single computation and \texttt{128/64} for a quad/double computation. All
experiments are performed on an \emph{Intel Xeon Platinum 8480+} CPU with
\emph{Julia  version 1.11.4}. The code supporting this work is openly available at \url{https://doi.org/10.5281/zenodo.22063645}.

\section{Performance Evaluation}\label{sec:results}
We evaluate the correction strategies on two nonlinear PDEs, semi-discretized in
space with $N_x$ Fourier spectral points and integrated in time with SDIRK2,
SDIRK3, and SDIRK4.

The two test problems considered are the inviscid Burgers' equation and the porous medium equation. For the inviscid Burgers' equation, defined by $u_t + \left(\tfrac12 u^2\right)_x = 0$ on $(0, 2\pi)$ with initial condition $u(x,0) = \sin(x)$ and periodic boundary conditions, the system is integrated to $T_f = 0.7$ (prior to shock formation) using the dominant operator $\mathcal{L} = -D_x$. For the porous medium equation, defined by $u_t = \left(u^3\right)_{xx}$ on $(-\pi, \pi)$ with $u(x,0) = \tfrac12\cos(x) + \tfrac{1}{2}$ and periodic boundary conditions, the system is integrated to $T_f = 0.5$ using the dominant operator $\mathcal{L} = D_{xx}$. In both cases, the errors in the $\ell^\infty$-norm are measured at the final time $T_f$.

\subsection{SDIRK2 (Implicit Midpoint Rule)}\label{sec:sdirk2}
We first study the second-order one-stage SDIRK2 (IMR) on Burgers' equation.
Table~\ref{tab:sdirk2} reports the final-time error at a fixed time step
$\dt = 10^{-3}$ for three spatial resolutions, $N_x \in \{50, 100, 200\}$, using
the precision pairings \texttt{64/16}, \texttt{64/32}, and \texttt{64/64}, for the
three correction strategies. The parenthesized values are the runtime speedups of
the mixed precision pairings over the full \texttt{64/64} computation for the same
$\dt$ and $N_x$. For this test case, the speedups of the mixed precision SDIRK2 improve as $N_x$ increases, 
regardless of the correction strategy employed. Furthermore, the \texttt{64/16} pairings are at least $8\times$ 
faster -- and up to $53\times$ faster -- than the implementations in full double precision.

\begin{table}[h]
  \caption{SDIRK2 on Burgers' equation ($\dt = 10^{-3}$): final-time error by spatial
  resolution and precision pairing, for each correction strategy. Parentheses
  give the speedup over \texttt{64/64}.}
  \label{tab:sdirk2}
  \centering
  \scriptsize
  \setlength{\tabcolsep}{3pt}
  \begin{tabular}{@{}c@{}}
    \begin{tabular}{@{}lccc@{}}
      \multicolumn{4}{c}{No corrections}\\
      \toprule
      $N_x$ & 64/16 & 64/32 & 64/64 \\
      \midrule
      50  & $9.06\times 10^{-4}$ ($24.3\times$) & $2.97 \times 10^{-7}$ ($1.9\times$) & $2.96 \times 10^{-7}$ \\
      100 & $8.82\times 10^{-4}$ ($38.0\times$) & $3.41 \times 10^{-7}$ ($1.9\times$) & $2.87 \times 10^{-7}$ \\
      200 & $8.83\times 10^{-4}$ ($72.7\times$) & $8.64 \times 10^{-7}$ ($2.1\times$) & $2.88 \times 10^{-7}$ \\
      \bottomrule
    \end{tabular}
    \\[1.5ex]
    \begin{tabular}{@{}lccc@{}}
      \multicolumn{4}{c}{One explicit correction}\\
      \toprule
      $N_x$ & 64/16 & 64/32 & 64/64 \\
      \midrule
      50  & $6.47 \times 10^{-7}$ ($10.1\times$) & $2.96 \times 10^{-7}$ ($2.2\times$) & $2.96 \times 10^{-7}$ \\
      100 & $8.44 \times 10^{-7}$ ($23.3\times$) & $2.87 \times 10^{-7}$ ($1.9\times$) & $2.87 \times 10^{-7}$ \\
      200 & $8.43 \times 10^{-7}$ ($63.6\times$) & $2.87 \times 10^{-7}$ ($2.1\times$) & $2.88 \times 10^{-7}$ \\
      \bottomrule
    \end{tabular}
    \\[1.5ex]
\begin{tabular}{@{}lccc@{}}
      \multicolumn{4}{c}{One $\Phi_{\mathrm{EIN}}$ correction}\\
      \toprule
      $N_x$ & 64/16 & 64/32 & 64/64 \\
      \midrule
      50  & $1.84 \times 10^{-7}$ ($8.9\times$)  & $2.96 \times 10^{-7}$ ($1.9\times$) & $2.96 \times 10^{-7}$ \\
      100 & $1.93 \times 10^{-6}$ ($20.7\times$) & $2.87 \times 10^{-7}$ ($1.9\times$) & $2.87 \times 10^{-7}$ \\
      200 & $2.20 \times 10^{-6}$ ($53.1\times$) & $2.87 \times 10^{-7}$ ($2.1\times$) & $2.88 \times 10^{-7}$ \\
      \bottomrule
    \end{tabular}
    \\[1.5ex]
    \begin{tabular}{@{}lccc@{}}
      \multicolumn{4}{c}{One $\Phi_J$ correction}\\
      \toprule
      $N_x$ & 64/16 & 64/32 & 64/64 \\
      \midrule
      50  & $2.39 \times 10^{-7}$ ($8.3\times$)  & $2.96 \times 10^{-7}$ ($1.8\times$) & $2.96 \times 10^{-7}$ \\
      100 & $2.24 \times 10^{-7}$ ($20.4\times$) & $2.87 \times 10^{-7}$ ($1.9\times$) & $2.87 \times 10^{-7}$ \\
      200 & $2.23 \times 10^{-7}$ ($53.9\times$) & $2.87 \times 10^{-7}$ ($2.1\times$) & $2.88 \times 10^{-7}$ \\
      \bottomrule
    \end{tabular}
  \end{tabular}
\end{table}




\subsection{SDIRK3}\label{sec:sdirk3}

For the third-order, two-stage SDIRK3 scheme applied to Burgers' equation, Table~\ref{tab:sdirk3nx} first reports the errors and speedups at $\dt= 10^{-3}$ and $N_x \in \{50, 100, 200\}$ for the precision pairings \texttt{64/16}, \texttt{64/32}, and \texttt{64/64}. As with the SDIRK2 scheme, we observe that the speedups improve for larger $N_x$ when combining double and half precision.

\begin{table}[t]
  \caption{SDIRK3 on Burgers' equation ($\dt = 10^{-3}$): final-time error by spatial
  resolution and precision pairing, for each correction strategy. Parentheses
  give the speedup over \texttt{64/64}.}
  \label{tab:sdirk3nx}
  \centering
  \scriptsize
  \setlength{\tabcolsep}{3pt}
  \begin{tabular}{@{}c@{}}
    \begin{tabular}{@{}lccc@{}}
      \multicolumn{4}{c}{No corrections}\\
      \toprule
      $N_x$ & 64/16 & 64/32 & 64/64 \\
      \midrule
      50  & $9.06 \times 10^{-4}$ ($14.7\times$) & $6.63 \times 10^{-9}$ ($2.0\times$) & $1.04 \times 10^{-9}$ \\
      100 & $8.82 \times 10^{-4}$ ($28.4\times$) & $2.36 \times 10^{-7}$ ($2.0\times$) & $1.44 \times 10^{-9}$ \\
      200 & $8.83 \times 10^{-4}$ ($94.0\times$) & $1.38 \times 10^{-6}$ ($2.1\times$) & $1.54 \times 10^{-9}$ \\
      \bottomrule
    \end{tabular}
    \\[1.5ex]
    \begin{tabular}{@{}lccc@{}}
      \multicolumn{4}{c}{Two explicit corrections}\\
      \toprule
      $N_x$ & 64/16 & 64/32 & 64/64 \\
      \midrule
      50  & $6.81 \times 10^{-9}$ ($11.8\times$) & $1.04 \times 10^{-9}$ ($1.8\times$) & $1.04 \times 10^{-9}$ \\
      100 & $7.58 \times 10^{-9}$ ($22.4\times$) & $1.44 \times 10^{-9}$ ($2.0\times$) & $1.44 \times 10^{-9}$ \\
      200 & $7.69 \times 10^{-9}$ ($70.4\times$) & $1.55 \times 10^{-9}$ ($2.1\times$) & $1.54 \times 10^{-9}$ \\
      \bottomrule
    \end{tabular}
    \\[1.5ex]
    \begin{tabular}{@{}lccc@{}}
      \multicolumn{4}{c}{Two $\Phi_{\mathrm{EIN}}$ corrections}\\
      \toprule
      $N_x$ & 64/16 & 64/32 & 64/64 \\
      \midrule
      50  & $2.33 \times 10^{-8}$ ($7.4\times$)  & $1.04 \times 10^{-9}$ ($1.8\times$) & $1.04 \times 10^{-9}$ \\
      100 & $2.76 \times 10^{-8}$ ($18.7\times$) & $1.45 \times 10^{-9}$ ($1.9\times$) & $1.44 \times 10^{-9}$ \\
      200 & $2.74 \times 10^{-8}$ ($49.6\times$) & $1.56 \times 10^{-9}$ ($2.0\times$) & $1.54 \times 10^{-9}$ \\
      \bottomrule
    \end{tabular}
    \\[1.5ex]
    \begin{tabular}{@{}lccc@{}}
      \multicolumn{4}{c}{Two $\Phi_J$ corrections}\\
      \toprule
      $N_x$ & 64/16 & 64/32 & 64/64 \\
      \midrule
      50  & $1.78 \times 10^{-9}$ ($7.4\times$)  & $1.04 \times 10^{-9}$ ($1.8\times$) & $1.04 \times 10^{-9}$ \\
      100 & $2.93 \times 10^{-9}$ ($19.9\times$) & $1.44 \times 10^{-9}$ ($2.0\times$) & $1.44 \times 10^{-9}$ \\
      200 & $3.41 \times 10^{-9}$ ($48.5\times$) & $1.54 \times 10^{-9}$ ($2.0\times$) & $1.54 \times 10^{-9}$ \\
      \bottomrule
    \end{tabular}
  \end{tabular}
\end{table}

Then we fix the spatial resolution and vary the time step, $\Delta t \in \{10^{-2}, 10^{-3}, 10^{-4}\}$, using the precision pairings \texttt{64/16}, \texttt{64/32}, and \texttt{64/64}. Table~\ref{tab:sdirk3} reports the final-time errors for the four correction strategies, with speedups relative to \texttt{64/64} for the same $N_x$ and $\Delta t$. While uncorrected mixed precision implementations lose some accuracy, the addition of any of the three correction strategies recovers this accuracy while maintaining significant runtime savings. Notably, the \texttt{64/16} pairing consistently achieves speedups of at least $7\times$.

\begin{table}[t]
  \caption{SDIRK3 on Burgers' equation ($N_x = 50$): final-time error by time step
  and precision pairing, for each correction strategy. Parentheses give the
  speedup over \texttt{64/64}.}
  \label{tab:sdirk3}
  \centering
  \scriptsize
  \setlength{\tabcolsep}{3pt}
  \begin{tabular}{@{}c@{}}
    \begin{tabular}{@{}lccc@{}}
      \multicolumn{4}{c}{No corrections}\\
      \toprule
     $\dt$ & 64/16 & 64/32 & 64/64 \\
      \midrule
      $10^{-2}$ & $8.90 \times 10^{-3}$ ($13.7\times$) & $1.64 \times 10^{-6}$ ($1.2\times$) & $9.69 \times 10^{-7}$ \\
      $10^{-3}$ & $9.06 \times 10^{-4}$ ($14.7\times$) & $6.63 \times 10^{-9}$ ($2.0\times$) & $1.04 \times 10^{-9}$ \\
      $10^{-4}$ & $9.08 \times 10^{-5}$ ($12.2\times$) & $7.45 \times 10^{-10}$ ($1.9\times$) & $1.09 \times 10^{-12}$ \\
      \bottomrule
    \end{tabular}
    \\[1.5ex]
    \begin{tabular}{@{}lccc@{}}
      \multicolumn{4}{c}{Two explicit corrections}\\
      \toprule
     $\dt$ & 64/16 & 64/32 & 64/64 \\
      \midrule
      $10^{-2}$ & $6.88 \times 10^{-6}$ ($10.6\times$)  & $9.68 \times 10^{-7}$ ($1.3\times$) & $9.69 \times 10^{-7}$ \\
      $10^{-3}$ & $6.81 \times 10^{-9}$ ($11.8\times$)  & $1.04 \times 10^{-9}$ ($1.8\times$) & $1.04 \times 10^{-9}$ \\
      $10^{-4}$ & $6.85 \times 10^{-12}$ ($10.0\times$) & $1.09 \times 10^{-12}$ ($1.8\times$) & $1.09 \times 10^{-12}$ \\
      \bottomrule
    \end{tabular}
    \\[1.5ex]
  \begin{tabular}{@{}lccc@{}}
      \multicolumn{4}{c}{Two $\Phi_{\mathrm{EIN}}$ corrections}\\
      \toprule
     $\dt$ & 64/16 & 64/32 & 64/64 \\
      \midrule
      $10^{-2}$ & $2.34 \times 10^{-5}$ ($8.4\times$)   & $9.69 \times 10^{-7}$ ($1.2\times$) & $9.69 \times 10^{-7}$ \\
      $10^{-3}$ & $2.33 \times 10^{-8}$ ($7.4\times$)   & $1.04 \times 10^{-9}$ ($1.8\times$) & $1.04 \times 10^{-9}$ \\
      $10^{-4}$ & $2.33 \times 10^{-11}$ ($7.8\times$)  & $1.09 \times 10^{-12}$ ($1.7\times$) & $1.09 \times 10^{-12}$ \\
      \bottomrule
    \end{tabular}
    \\[1.5ex]
    \begin{tabular}{@{}lccc@{}}
      \multicolumn{4}{c}{Two $\Phi_J$ corrections}\\
      \toprule
     $\dt$ & 64/16 & 64/32 & 64/64 \\
      \midrule
      $10^{-2}$ & $1.82 \times 10^{-6}$ ($8.3\times$)   & $9.69 \times 10^{-7}$ ($1.2\times$) & $9.69 \times 10^{-7}$ \\
      $10^{-3}$ & $1.78 \times 10^{-9}$ ($7.4\times$)   & $1.04 \times 10^{-9}$ ($1.8\times$) & $1.04 \times 10^{-9}$ \\
      $10^{-4}$ & $1.82 \times 10^{-12}$ ($7.8\times$)  & $1.09 \times 10^{-12}$ ($1.7\times$) & $1.09 \times 10^{-12}$ \\
      \bottomrule
    \end{tabular}
  \end{tabular}
\end{table}

\begin{figure}[t]
  \centering
  \begin{subfigure}{0.48\columnwidth}
    \centering
    \includegraphics[width=\linewidth]{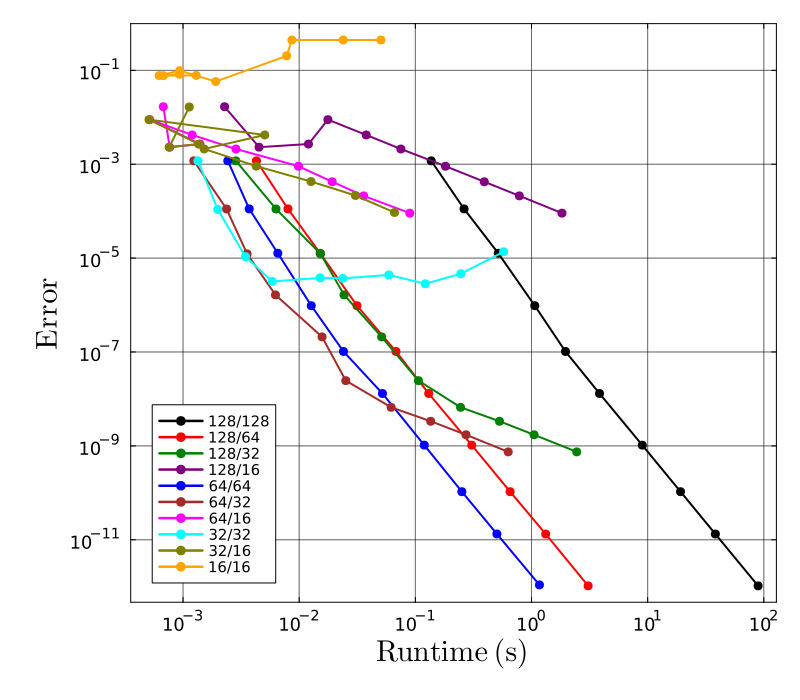}
  \end{subfigure}\hfill
  \begin{subfigure}{0.48\columnwidth}
    \centering
    \includegraphics[width=\linewidth]{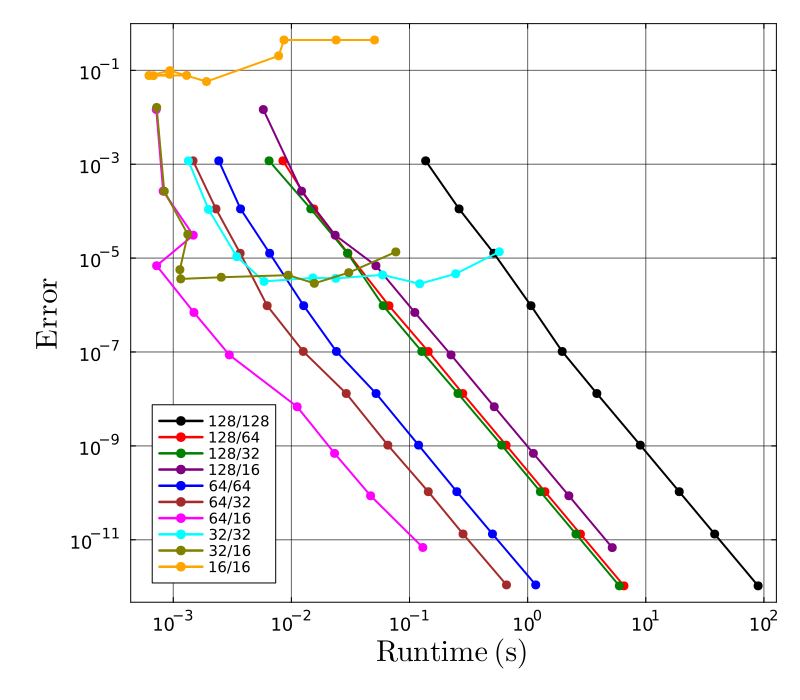}
  \end{subfigure}

  \vspace{1ex}

\begin{subfigure}{0.48\columnwidth}
    \centering
    \includegraphics[width=\linewidth]{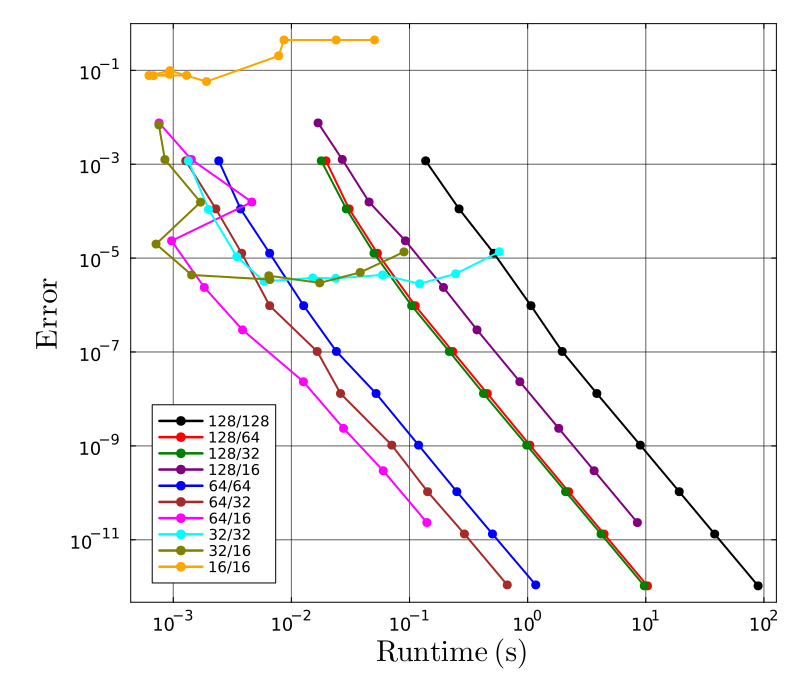}
  \end{subfigure}\hfill
  \begin{subfigure}{0.48\columnwidth}
    \centering
    \includegraphics[width=\linewidth]{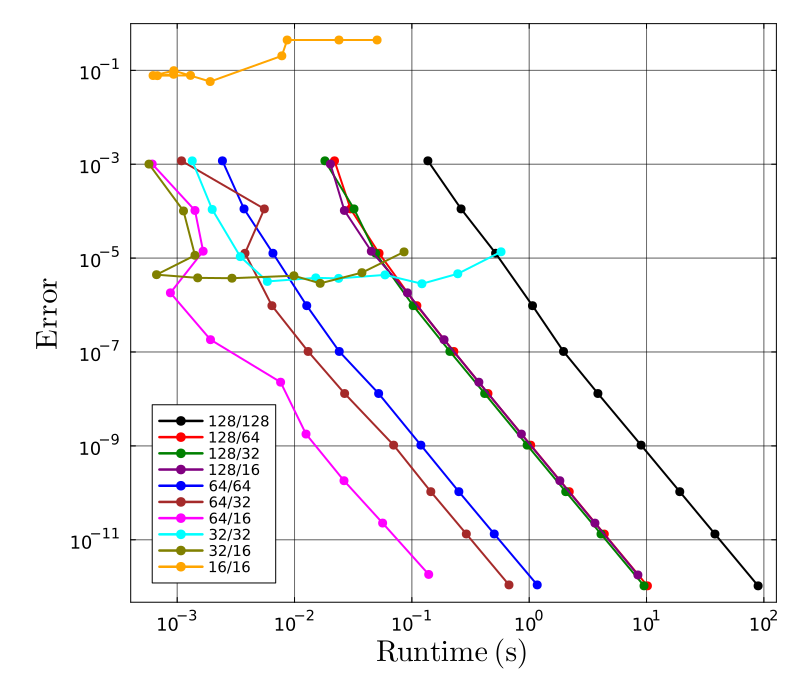}
  \end{subfigure}
  \caption{SDIRK3 on Burgers' equation: error vs.\ runtime for each correction strategy.
  Top left: no corrections. Top right: two explicit corrections. Bottom left: two
  $\Phi_{\mathrm{EIN}}$ corrections. Bottom right: two $\Phi_J$ corrections.}
  \label{fig:sdirk3}
\end{figure}

Figure~\ref{fig:sdirk3} shows the error against runtime for SDIRK3 on Burgers' equation under the three correction strategies for $N_x = 50$. Similarly to the SDIRK2 method, combining double and half precision with the application of $\Phi_J$ corrections recovers accuracy close to that of the full-precision implementations while achieving significant runtime savings.

\subsection{SDIRK4}

\begin{figure}[t]
  \centering
  \begin{subfigure}{0.48\columnwidth}
    \centering
    \includegraphics[width=\linewidth]{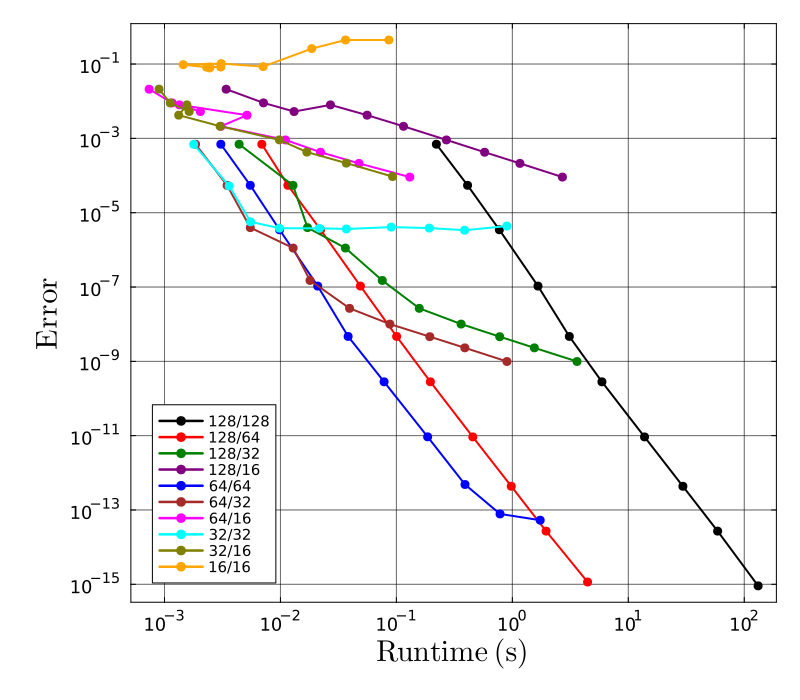}
  \end{subfigure}\hfill
  \begin{subfigure}{0.48\columnwidth}
    \centering
    \includegraphics[width=\linewidth]{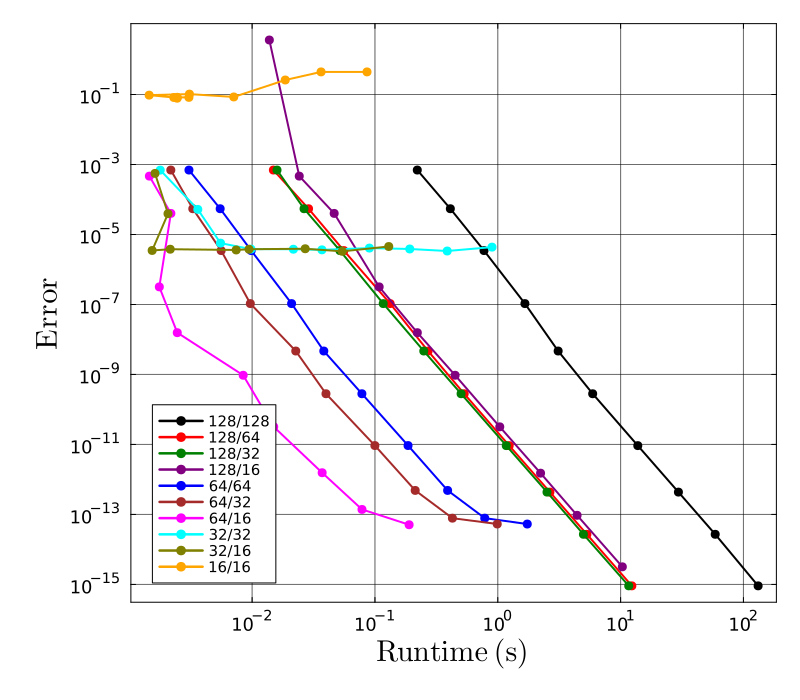}
  \end{subfigure}

  \vspace{1ex}

\begin{subfigure}{0.48\columnwidth}
    \centering
    \includegraphics[width=\linewidth]{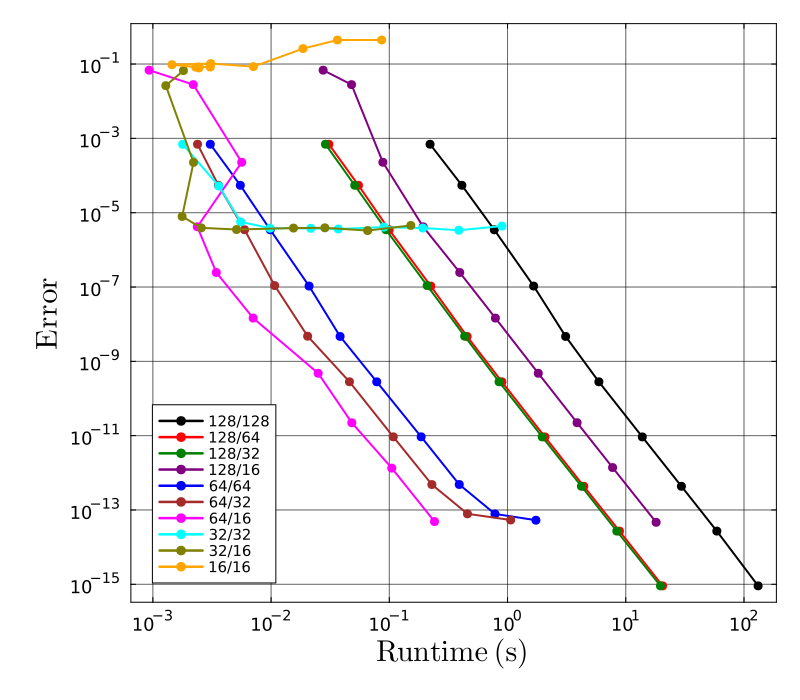}
  \end{subfigure}\hfill
  \begin{subfigure}{0.48\columnwidth}
    \centering
    \includegraphics[width=\linewidth]{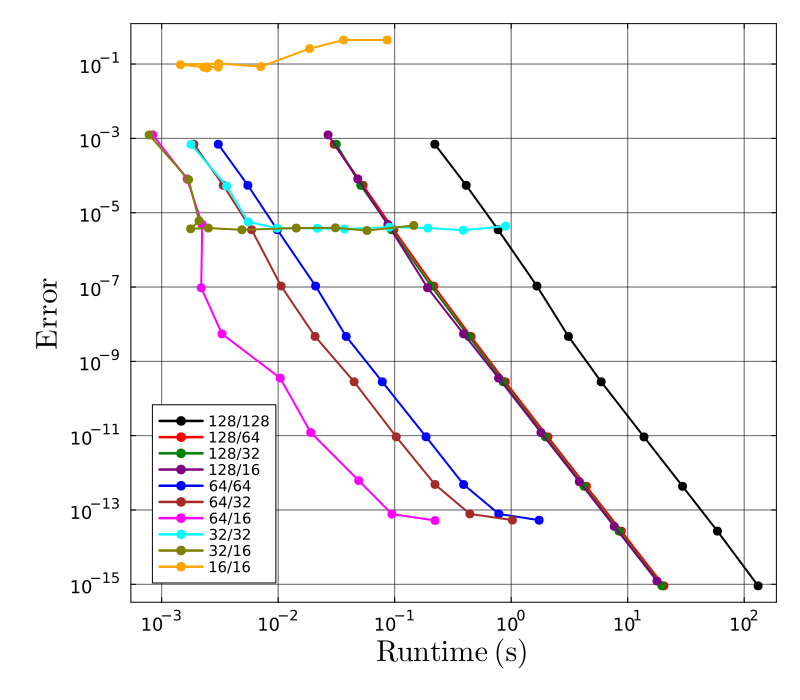}
  \end{subfigure}
  \caption{SDIRK4 on Burgers' equation: error vs.\ runtime for
  each correction strategy. Top left: no corrections. Top right: three explicit corrections. Bottom left: three
  $\Phi_{\mathrm{EIN}}$ corrections. Bottom right: three $\Phi_J$ corrections.}
  \label{fig:sdirk4Burgers}
\end{figure}

For the fourth-order, three-stage SDIRK4 method applied to Burgers' equation, we again fix the spatial resolution $N_x$ and compare several quadruple-precision variations across various time steps, $\dt$. Table~\ref{tab:sdirk4quadBurgers} compares the mixed precision pairings \texttt{128/16} and \texttt{128/32} against the full-precision \texttt{128/128} implementation. Figure~\ref{fig:sdirk4Burgers} plots runtime versus error for a wider variety of precision combinations across the different correction approaches. The results demonstrate that the $\Phi_J$ correction strategy recovers the most accuracy for the mixed precision implementations.

\begin{table}[t]
  \caption{SDIRK4 on Burgers' equation ($N_x = 50$): final-time error
  by time step and precision pairing, for each correction strategy. Parentheses
  give the speedup over \texttt{128/128}.}
  \label{tab:sdirk4quadBurgers}
  \centering
  \scriptsize
  \setlength{\tabcolsep}{3pt}
  \begin{tabular}{@{}c@{}}
    \begin{tabular}{@{}lccc@{}}
      \multicolumn{4}{c}{No corrections}\\
      \toprule
     $\dt$ & 128/16 & 128/32 & 128/128 \\
      \midrule
      $10^{-2}$ & $7.96 \times 10^{-3}$ ($56.6\times$) & $1.12 \times 10^{-6}$ ($38.9\times$) & $1.06 \times 10^{-7}$ \\
      $10^{-3}$ & $9.06 \times 10^{-4}$ ($48.7\times$) & $1.01 \times 10^{-8}$ ($35.7\times$) & $9.26 \times 10^{-12}$ \\
      $10^{-4}$ & $9.08 \times 10^{-5}$ ($46.3\times$) & $9.90 \times 10^{-10}$ ($34.9\times$) & $9.10 \times 10^{-16}$ \\
      \bottomrule
    \end{tabular}
    \\[1.5ex]
    \begin{tabular}{@{}lccc@{}}
      \multicolumn{4}{c}{Three explicit corrections}\\
      \toprule
     $\dt$ & 128/16 & 128/32 & 128/128 \\
      \midrule
      $10^{-2}$ & $3.18 \times 10^{-7}$ ($15.2\times$)  & $1.06 \times 10^{-7}$ ($13.7\times$) & $1.06 \times 10^{-7}$ \\
      $10^{-3}$ & $3.19 \times 10^{-11}$ ($13.2\times$) & $9.26 \times 10^{-12}$ ($11.7\times$) & $9.26 \times 10^{-12}$ \\
      $10^{-4}$ & $3.18 \times 10^{-15}$ ($12.7\times$) & $9.10 \times 10^{-16}$ ($11.3\times$) & $9.10 \times 10^{-16}$ \\
      \bottomrule
    \end{tabular}
    \\[1.5ex]
  \begin{tabular}{@{}lccc@{}}
      \multicolumn{4}{c}{Three $\Phi_{\mathrm{EIN}}$ corrections}\\
      \toprule
      $\dt$ & 128/16 & 128/32 & 128/128 \\
      \midrule
      $10^{-2}$ & $4.17 \times 10^{-6}$ ($8.5\times$)  & $1.09 \times 10^{-7}$ ($7.8\times$) & $1.06 \times 10^{-7}$ \\
      $10^{-3}$ & $4.77 \times 10^{-10}$ ($7.4\times$) & $9.26 \times 10^{-12}$ ($6.9\times$) & $9.26 \times 10^{-12}$ \\
      $10^{-4}$ & $4.66 \times 10^{-14}$ ($7.2\times$) & $9.10 \times 10^{-16}$ ($6.6\times$) & $9.10 \times 10^{-16}$ \\
      \bottomrule
    \end{tabular}
    \\[1.5ex]
    \begin{tabular}{@{}lccc@{}}
      \multicolumn{4}{c}{Three $\Phi_J$ corrections}\\
      \toprule
     $\dt$ & 128/16 & 128/32 & 128/128 \\
      \midrule
      $10^{-2}$ & $9.53 \times 10^{-8}$ ($8.8\times$)  & $1.06 \times 10^{-7}$ ($7.8\times$) & $1.06 \times 10^{-7}$ \\
      $10^{-3}$ & $1.21 \times 10^{-11}$ ($7.6\times$) & $9.26 \times 10^{-12}$ ($6.8\times$) & $9.26 \times 10^{-12}$ \\
      $10^{-4}$ & $1.23 \times 10^{-15}$ ($7.4\times$) & $9.10 \times 10^{-16}$ ($6.7\times$) & $9.10 \times 10^{-16}$ \\
      \bottomrule
    \end{tabular}
  \end{tabular}
\end{table}

\begin{figure}[t]
  \centering
  \begin{subfigure}{0.48\columnwidth}
    \centering
    \includegraphics[width=\linewidth]{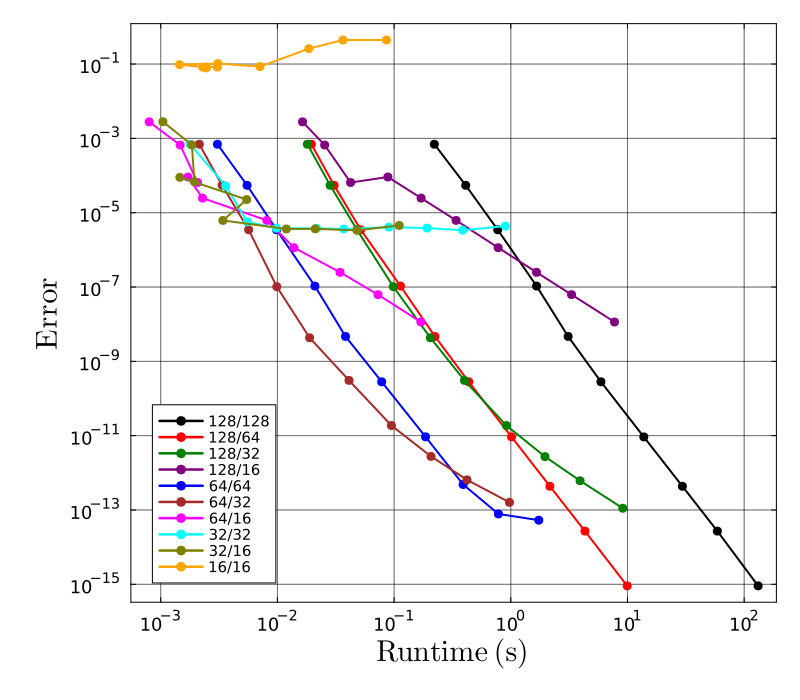}
  \end{subfigure}\hfill
  \begin{subfigure}{0.48\columnwidth}
    \centering
    \includegraphics[width=\linewidth]{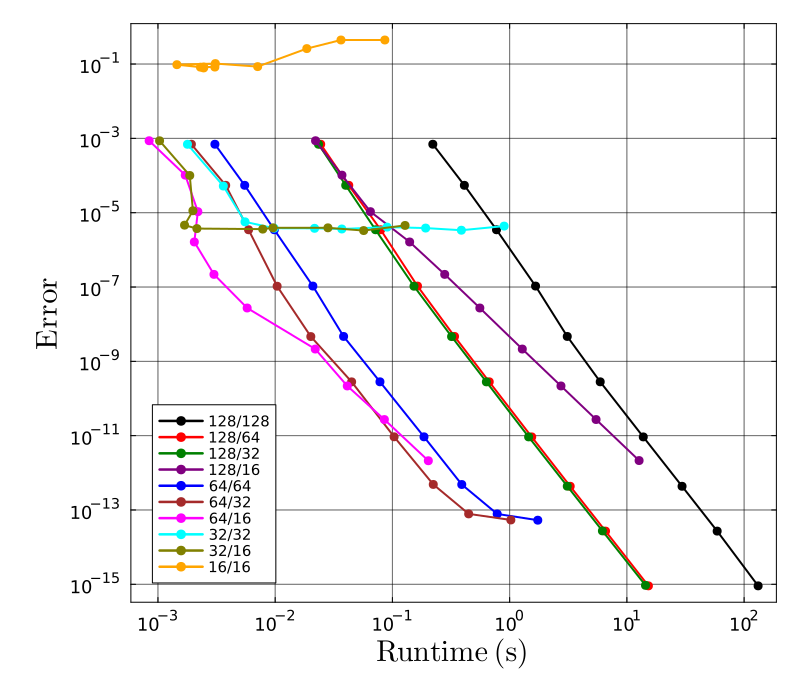}
  \end{subfigure}

  \vspace{1ex}

  \begin{subfigure}{0.48\columnwidth}
    \centering
    \includegraphics[width=\linewidth]{Figures/Burgers/Figure_PlotTime_SDIRK4_Nx50_3_JAC.png}
  \end{subfigure}\hfill
  \begin{subfigure}{0.48\columnwidth}
    \centering
    \includegraphics[width=\linewidth]{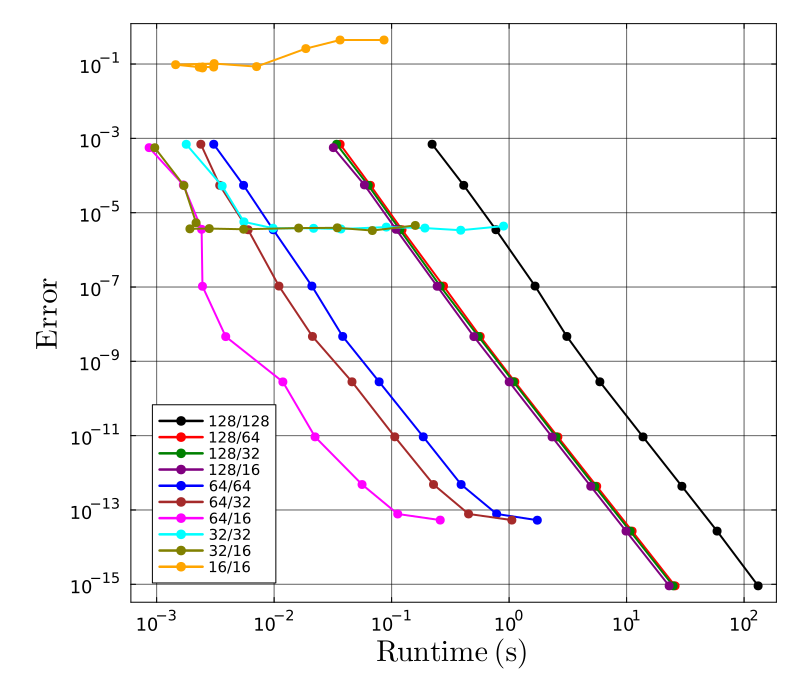}
  \end{subfigure}
  \caption{SDIRK4 on Burgers' equation: error vs.\ runtime for
  increasing numbers of $\Phi_J$ corrections. Top left: one $\Phi_J$ correction. Top right:
  two $\Phi_J$ corrections. Bottom left: three $\Phi_J$ corrections. Bottom right:
  four $\Phi_J$ corrections.}
  \label{fig:sdirk4B}
\end{figure}

We evaluate the performance when varying the number of $\Phi_J$ corrections. Table~\ref{tab:sdirk4phijBurgers} shows the errors and speedups for one to four corrections using the mixed precision combinations \texttt{128/16} and \texttt{128/32}. We observe that the error improves remarkably as the number of corrections increases from one to four, though this extra accuracy comes at the cost of additional runtime. However, speedups of more than $5\times$ over the \texttt{128/128} implementation are still obtained even when applying four corrections. The runtime versus error for all precision combinations across these different numbers of $\Phi_J$ corrections is presented in Figure~\ref{fig:sdirk4B}.

\begin{table}[t]
  \caption{SDIRK4 on Burgers' equation ($N_x = 50$): final-time
  error by time step and precision pairing, for increasing numbers of $\Phi_J$
  corrections. Parentheses give the speedup over \texttt{128/128}.}
  \label{tab:sdirk4phijBurgers}
  \centering
  \scriptsize
  \setlength{\tabcolsep}{3pt}
  \begin{tabular}{@{}c@{}}
    \begin{tabular}{@{}lccc@{}}
      \multicolumn{4}{c}{One $\Phi_J$ corr.}\\
      \toprule
      $\dt$ & 128/16 & 128/32 & 128/128 \\
      \midrule
      $10^{-2}$ & $9.11 \times 10^{-5}$ ($17.9\times$)   & $1.02 \times 10^{-7}$ ($15.6\times$) & $1.06 \times 10^{-7}$ \\
      $10^{-3}$ & $1.14 \times 10^{-6}$ ($17.0\times$)   & $1.87 \times 10^{-11}$ ($14.5\times$) & $9.26 \times 10^{-12}$ \\
      $10^{-4}$ & $1.14 \times 10^{-8}$ ($16.6\times$)   & $1.11 \times 10^{-13}$ ($14.1\times$) & $9.10 \times 10^{-16}$ \\
      \bottomrule
    \end{tabular}
    \\[1.5ex]
    \begin{tabular}{@{}lccc@{}}
      \multicolumn{4}{c}{Two $\Phi_J$ corr.}\\
      \toprule
      $\dt$ & 128/16 & 128/32 & 128/128 \\
      \midrule
      $10^{-2}$ & $1.63 \times 10^{-6}$ ($11.4\times$)   & $1.05 \times 10^{-7}$ ($10.2\times$) & $1.06 \times 10^{-7}$ \\
      $10^{-3}$ & $2.14 \times 10^{-9}$ ($10.3\times$)   & $9.27 \times 10^{-12}$ ($9.1\times$)  & $9.26 \times 10^{-12}$ \\
      $10^{-4}$ & $2.14 \times 10^{-12}$ ($10.0\times$)  & $9.39 \times 10^{-16}$ ($8.8\times$)  & $9.10 \times 10^{-16}$ \\
      \bottomrule
    \end{tabular}
    \\[1.5ex]
    \begin{tabular}{@{}lccc@{}}
      \multicolumn{4}{c}{Three $\Phi_J$ corr.}\\
      \toprule
      $\dt$ & 128/16 & 128/32 & 128/128 \\
      \midrule
      $10^{-2}$ & $9.53 \times 10^{-8}$ ($8.8\times$)   & $1.06 \times 10^{-7}$ ($7.8\times$) & $1.06 \times 10^{-7}$ \\
      $10^{-3}$ & $1.21 \times 10^{-11}$ ($7.6\times$)  & $9.26 \times 10^{-12}$ ($6.8\times$) & $9.26 \times 10^{-12}$ \\
      $10^{-4}$ & $1.23 \times 10^{-15}$ ($7.4\times$)  & $9.10 \times 10^{-16}$ ($6.7\times$) & $9.10 \times 10^{-16}$ \\
      \bottomrule
    \end{tabular}
    \\[1.5ex]
    \begin{tabular}{@{}lccc@{}}
      \multicolumn{4}{c}{Four $\Phi_J$ corr.}\\
      \toprule
      $\dt$ & 128/16 & 128/32 & 128/128 \\
      \midrule
      $10^{-2}$ & $1.05 \times 10^{-7}$ ($6.6\times$)   & $1.06 \times 10^{-7}$ ($6.1\times$) & $1.06 \times 10^{-7}$ \\
      $10^{-3}$ & $9.24 \times 10^{-12}$ ($5.7\times$)  & $9.26 \times 10^{-12}$ ($5.3\times$) & $9.26 \times 10^{-12}$ \\
      $10^{-4}$ & $9.10 \times 10^{-16}$ ($5.5\times$)  & $9.10 \times 10^{-16}$ ($5.1\times$) & $9.10 \times 10^{-16}$ \\
      \bottomrule
    \end{tabular}
  \end{tabular}
\end{table}

\begin{figure}[t]
  \centering
  \begin{subfigure}{0.48\columnwidth}
    \centering
    \includegraphics[width=\linewidth]{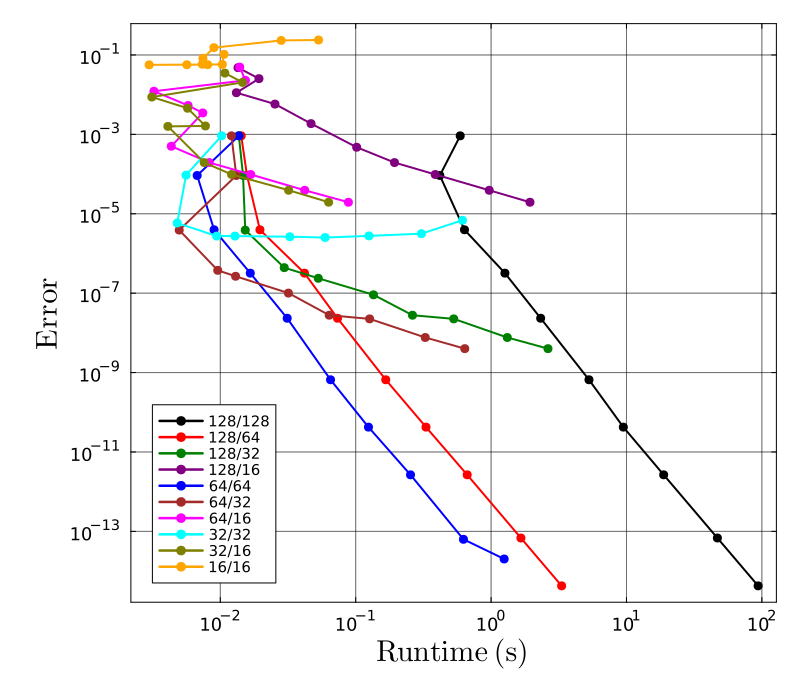}
  \end{subfigure}\hfill
  \begin{subfigure}{0.48\columnwidth}
    \centering
    \includegraphics[width=\linewidth]{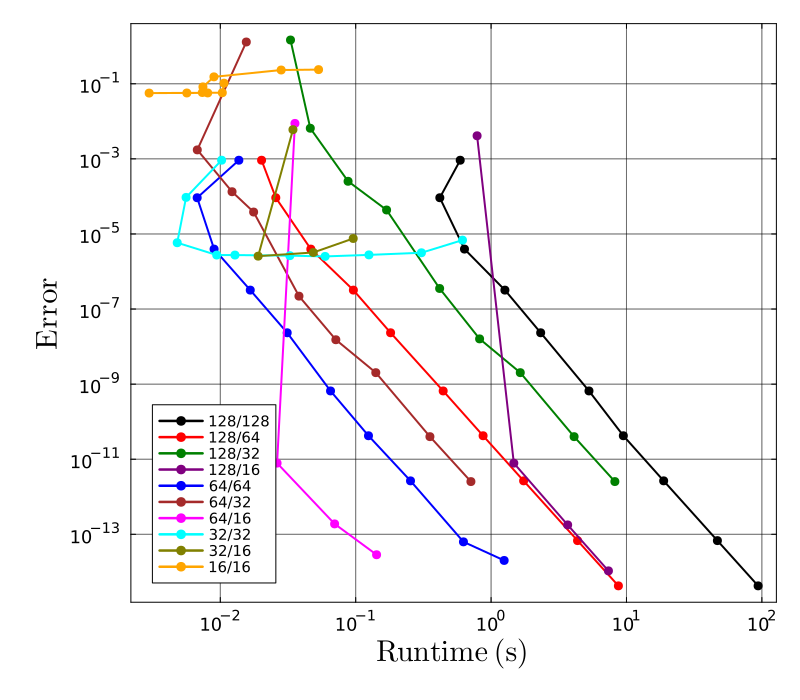}
  \end{subfigure}

  \vspace{1ex}

\begin{subfigure}{0.48\columnwidth}
    \centering
    \includegraphics[width=\linewidth]{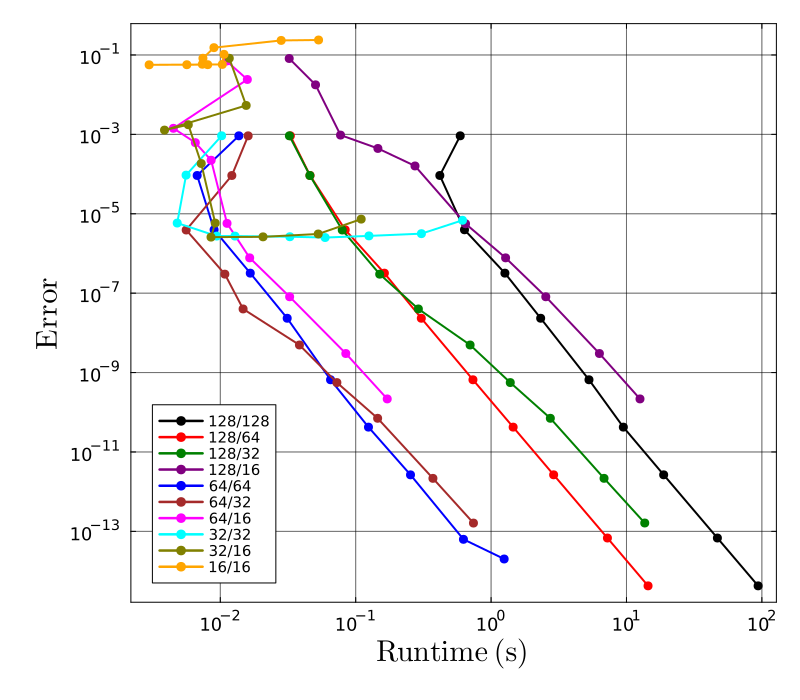}
  \end{subfigure}\hfill
  \begin{subfigure}{0.48\columnwidth}
    \centering
    \includegraphics[width=\linewidth]{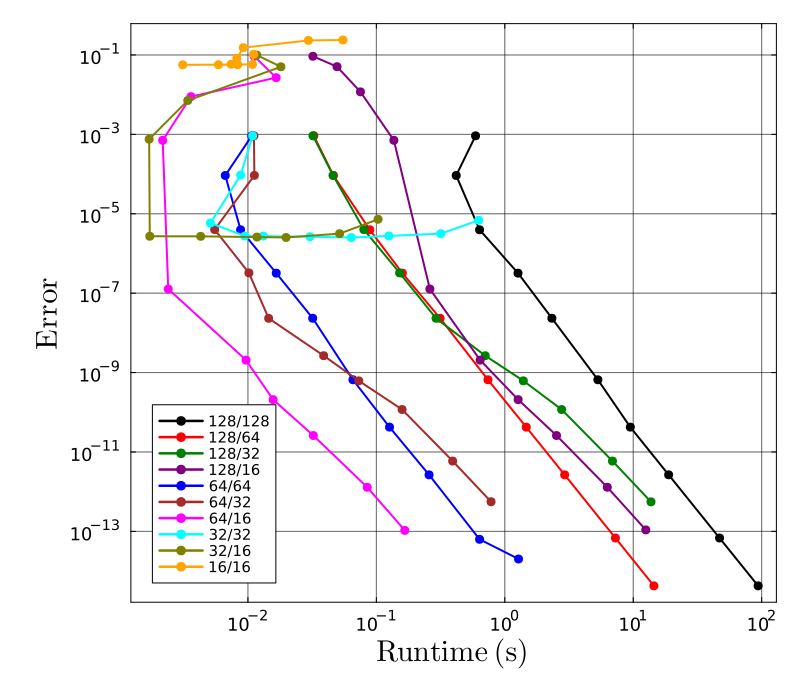}
  \end{subfigure}
  \caption{SDIRK4 on the porous medium equation: error vs.\ runtime for each correction strategy. Top left: no corrections. Top right: three explicit corrections. Bottom left: three
  $\Phi_{\mathrm{EIN}}$ corrections. Bottom right: three $\Phi_J$ corrections.}
  \label{fig:sdirk4}
\end{figure}

We conclude our performance evaluation by testing the mixed precision SDIRK4 method on the porous medium equation. Similar to the results for Burgers' equation, Table~\ref{tab:sdirk4} compares the mixed precision pairings \texttt{128/16} and \texttt{128/32} against the full-precision \texttt{128/128} implementation. For this problem, explicit corrections become unstable for larger $\dt$ when using half precision. In contrast, the $\Phi_{\mathrm{EIN}}$ and $\Phi_{J}$ corrections remain stable for all considered time steps $\dt$.  Figure~\ref{fig:sdirk4} illustrates runtime versus errors of other precision combinations for the different correction approaches. 

\begin{table}[t]
  \caption{SDIRK4 on the porous medium equation ($N_x = 50$): final-time error
  by time step and precision pairing, for each correction strategy. Parentheses
  give the speedup over \texttt{128/128}.}
  \label{tab:sdirk4}
  \centering
  \scriptsize
  \setlength{\tabcolsep}{3pt}
  \begin{tabular}{@{}c@{}}
    \begin{tabular}{@{}lccc@{}}
      \multicolumn{4}{c}{No corrections}\\
      \toprule
      $\dt$ & 128/16 & 128/32 & 128/128 \\
      \midrule
      $10^{-2}$ & $5.79 \times 10^{-3}$ ($50.1\times$) & $4.39 \times 10^{-7}$ ($43.7\times$) & $3.19 \times 10^{-7}$ \\
      $10^{-3}$ & $1.95 \times 10^{-4}$ ($49.0\times$) & $2.78 \times 10^{-8}$ ($37.0\times$) & $4.22 \times 10^{-11}$ \\
      $10^{-4}$ & $1.95 \times 10^{-5}$ ($48.9\times$) & $4.00 \times 10^{-9}$ ($36.4\times$) & $4.24 \times 10^{-15}$ \\
      \bottomrule
    \end{tabular}
    \\[1.5ex]
    \begin{tabular}{@{}lccc@{}}
      \multicolumn{4}{c}{Three explicit corrections}\\
      \toprule
      $\dt$ & 128/16 & 128/32 & 128/128 \\
      \midrule
      $10^{-2}$ & N/A                    & $2.53 \times 10^{-4}$ ($13.8\times$) & $3.19 \times 10^{-7}$ \\
      $10^{-3}$ & $4.12 \times 10^{-3}$ ($12.1\times$) & $1.60 \times 10^{-8}$ ($11.6\times$) & $4.22 \times 10^{-11}$ \\
      $10^{-4}$ & $1.06 \times 10^{-14}$ ($13.0\times$) & $2.54 \times 10^{-12}$ ($11.6\times$) & $4.24 \times 10^{-15}$ \\
      \bottomrule
    \end{tabular}
    \\[1.5ex]
  \begin{tabular}{@{}lccc@{}}
      \multicolumn{4}{c}{Three $\Phi_{\mathrm{EIN}}$ corrections}\\
      \toprule
      $\dt$ & 128/16 & 128/32 & 128/128 \\
      \midrule
      $10^{-2}$ & $4.42 \times 10^{-4}$ ($8.5\times$)  & $3.01 \times 10^{-7}$ ($8.4\times$) & $3.19 \times 10^{-7}$ \\
      $10^{-3}$ & $7.76 \times 10^{-7}$ ($7.5\times$)  & $5.60 \times 10^{-10}$ ($6.8\times$) & $4.22 \times 10^{-11}$ \\
      $10^{-4}$ & $2.17 \times 10^{-10}$ ($7.5\times$) & $1.62 \times 10^{-13}$ ($6.9\times$) & $4.24 \times 10^{-15}$ \\
      \bottomrule
    \end{tabular}
    \\[1.5ex]
    \begin{tabular}{@{}lccc@{}}
      \multicolumn{4}{c}{Three $\Phi_J$ corrections}\\
      \toprule
      $\dt$ & 128/16 & 128/32 & 128/128 \\
      \midrule
      $10^{-2}$ & $7.09 \times 10^{-4}$ ($8.9\times$)  & $3.24 \times 10^{-7}$ ($8.3\times$) & $3.19 \times 10^{-7}$ \\
      $10^{-3}$ & $2.08 \times 10^{-10}$ ($7.5\times$) & $6.19 \times 10^{-10}$ ($6.8\times$) & $4.22 \times 10^{-11}$ \\
      $10^{-4}$ & $1.09 \times 10^{-13}$ ($7.5\times$) & $5.55 \times 10^{-13}$ ($6.8\times$) & $4.24 \times 10^{-15}$ \\
      \bottomrule
    \end{tabular}
  \end{tabular}
\end{table}

\begin{figure}[t]
  \centering
  \begin{subfigure}{0.48\columnwidth}
    \centering
    \includegraphics[width=\linewidth]{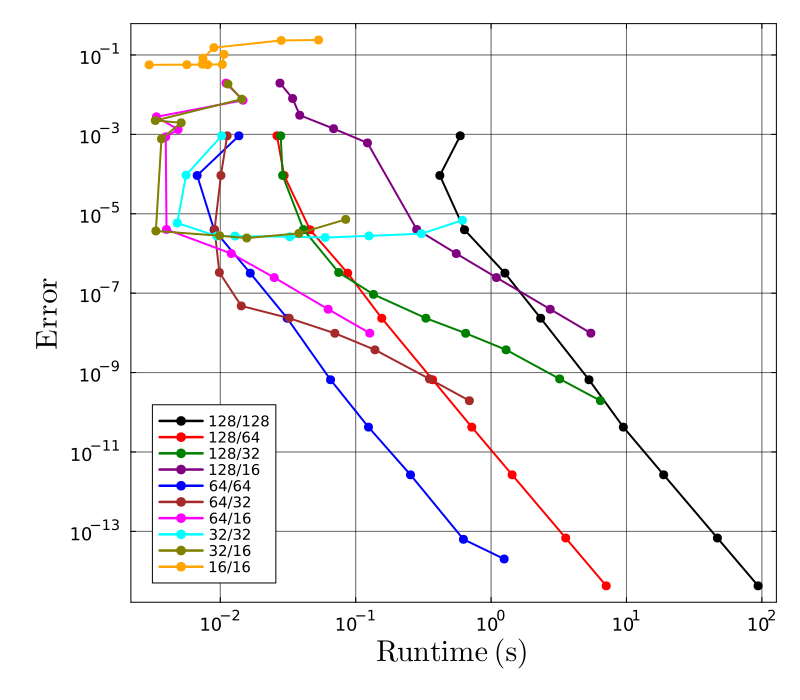}
  \end{subfigure}\hfill
  \begin{subfigure}{0.48\columnwidth}
    \centering
    \includegraphics[width=\linewidth]{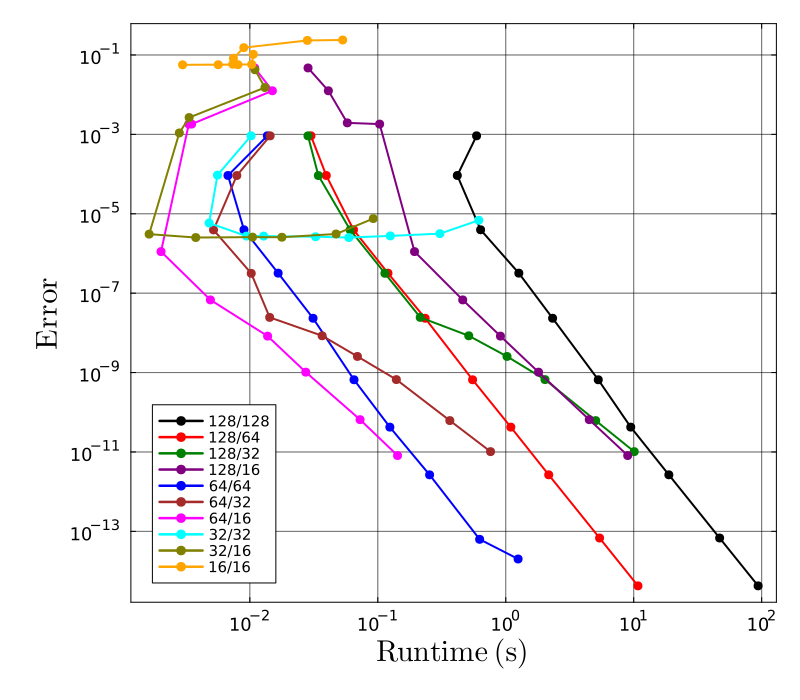}
  \end{subfigure}

  \vspace{1ex}

  \begin{subfigure}{0.48\columnwidth}
    \centering
    \includegraphics[width=\linewidth]{Figures/PM/Figure_PlotTime_SDIRK4_Nx50_3_JAC.png}
  \end{subfigure}\hfill
  \begin{subfigure}{0.48\columnwidth}
    \centering
    \includegraphics[width=\linewidth]{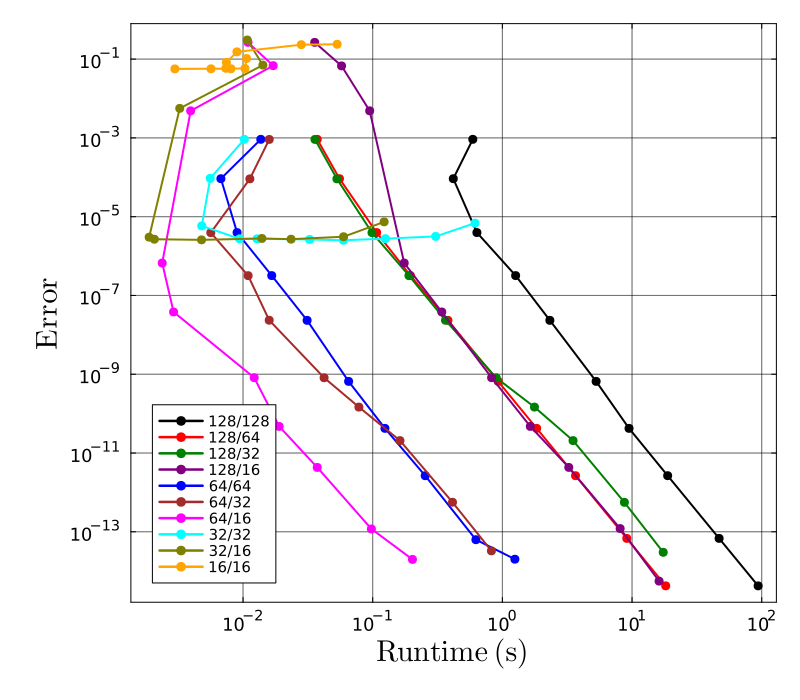}
  \end{subfigure}
  \caption{SDIRK4 on the porous medium equation: error vs.\ runtime for
  increasing numbers of $\Phi_J$ corrections. Top left: one $\Phi_J$ correction. Top right:
  two $\Phi_J$ corrections. Bottom left: three $\Phi_J$ corrections. Bottom right:
  four $\Phi_J$ corrections.}
  \label{fig:sdirk42}
\end{figure}

Finally, we evaluate the performance when varying the number of $\Phi_J$ corrections. As before, Table~\ref{tab:sdirk4phij} shows a remarkable improvement in accuracy from one to four corrections, with speedups comparable to those obtained for Burgers' equation. The runtime versus error for all precision combinations across these correction counts is presented in Figure~\ref{fig:sdirk42}.

\begin{table}[t]
  \caption{SDIRK4 on the porous medium equation ($N_x = 50$): final-time
  error by time step and precision pairing, for increasing numbers of $\Phi_J$
  corrections. Parentheses give the speedup over \texttt{128/128}.}
  \label{tab:sdirk4phij}
  \centering
  \scriptsize
  \setlength{\tabcolsep}{3pt}
  \begin{tabular}{@{}c@{}}
    \begin{tabular}{@{}lccc@{}}
      \multicolumn{4}{c}{One $\Phi_J$ corr.}\\
      \toprule
      $\dt$ & 128/16 & 128/32 & 128/128 \\
      \midrule
      $10^{-2}$ & $1.39 \times 10^{-3}$ ($18.4\times$)  & $3.33 \times 10^{-7}$ ($16.6\times$) & $3.19 \times 10^{-7}$ \\
      $10^{-3}$ & $9.94 \times 10^{-7}$ ($17.0\times$)  & $9.76 \times 10^{-9}$ ($14.6\times$) & $4.22 \times 10^{-11}$ \\
      $10^{-4}$ & $9.84 \times 10^{-9}$ ($17.1\times$)  & $1.96 \times 10^{-10}$ ($14.7\times$) & $4.24 \times 10^{-15}$ \\
      \bottomrule
    \end{tabular}
    \\[1.5ex]
    \begin{tabular}{@{}lccc@{}}
      \multicolumn{4}{c}{Two $\Phi_J$ corr.}\\
      \toprule
      $\dt$ & 128/16 & 128/32 & 128/128 \\
      \midrule
      $10^{-2}$ & $1.82 \times 10^{-3}$ ($12.3\times$)  & $3.17 \times 10^{-7}$ ($10.8\times$) & $3.19 \times 10^{-7}$ \\
      $10^{-3}$ & $8.29 \times 10^{-9}$ ($10.5\times$)  & $2.55 \times 10^{-9}$ ($9.3\times$)  & $4.22 \times 10^{-11}$ \\
      $10^{-4}$ & $8.13 \times 10^{-12}$ ($10.5\times$) & $1.02 \times 10^{-11}$ ($9.4\times$) & $4.24 \times 10^{-15}$ \\
      \bottomrule
    \end{tabular}
    \\[1.5ex]
    \begin{tabular}{@{}lccc@{}}
      \multicolumn{4}{c}{Three $\Phi_J$ corr.}\\
      \toprule
      $\dt$ & 128/16 & 128/32 & 128/128 \\
      \midrule
      $10^{-2}$ & $7.09 \times 10^{-4}$ ($9.2\times$)   & $3.24 \times 10^{-7}$ ($8.4\times$) & $3.19 \times 10^{-7}$ \\
      $10^{-3}$ & $2.08 \times 10^{-10}$ ($7.5\times$)  & $6.19 \times 10^{-10}$ ($6.8\times$) & $4.22 \times 10^{-11}$ \\
      $10^{-4}$ & $1.09 \times 10^{-13}$ ($7.5\times$)  & $5.55 \times 10^{-13}$ ($6.9\times$) & $4.24 \times 10^{-15}$ \\
      \bottomrule
    \end{tabular}
    \\[1.5ex]
    \begin{tabular}{@{}lccc@{}}
      \multicolumn{4}{c}{Four $\Phi_J$ corr.}\\
      \toprule
      $\dt$ & 128/16 & 128/32 & 128/128 \\
      \midrule
      $10^{-2}$ & $6.67 \times 10^{-7}$ ($7.2\times$)   & $3.18 \times 10^{-7}$ ($6.5\times$) & $3.19 \times 10^{-7}$ \\
      $10^{-3}$ & $4.73 \times 10^{-11}$ ($5.8\times$)  & $1.47 \times 10^{-10}$ ($5.4\times$) & $4.22 \times 10^{-11}$ \\
      $10^{-4}$ & $5.60 \times 10^{-15}$ ($5.8\times$)  & $3.01 \times 10^{-14}$ ($5.4\times$) & $4.24 \times 10^{-15}$ \\
      \bottomrule
    \end{tabular}
  \end{tabular}
\end{table}

\section{Interpretation of the observed speedups}
The large speedups observed for the half-precision pairings are consistent with the precision-dependent compute roof of the Intel Xeon Platinum 8480+ platform, but they should not be interpreted as a simple factor-of-four reduction in storage cost.  The 8480+ has 56 cores, a 2.00~GHz base frequency, and two AVX-512 fused-multiply-add (FMA) units per core \cite{Intel8480Spec}; moreover, 4th generation Intel Xeon Scalable processors provide native AVX-512 FP16 support for IEEE half-precision arithmetic \cite{IntelAVX512FP16}.  For a vector FMA kernel, the ideal peak rate can be estimated as
\[
P_{\max}
=
N_{\rm core} f_{Chip\, Clock} N_{\rm FMA}
\left(\frac{512}{b}\right)2,
\]
where \(b\) is the number of bits per floating-point value and the final factor of two counts the multiply and add in an FMA (fuze multiply add).  Thus a 512-bit AVX-512 register contains \(512/64=8\) double-precision values but \(512/16=32\) half-precision values.  At the 2.00~GHz base clock, this gives an ideal FP64 peak of approximately
\[
56 \cdot 2.00{\times}10^9 \cdot 2 \cdot 8 \cdot 2
\approx 3.6~{\rm TFLOP/s},
\]
whereas the corresponding IEEE FP16 AVX-512 peak is approximately
\[
56 \cdot 2.00{\times}10^9 \cdot 2 \cdot 32 \cdot 2
\approx 14.3~{\rm TFLOP/s}.
\]
This is a SIMD register-lane effect: more operands fit in each fixed-width vector register and therefore more arithmetic operations can be issued per vector instruction.  In addition, half precision reduces data movement and improves cache residency because a half-precision value occupies two bytes rather than eight.  The speedups exceeding the nominal AVX-512 FP16/FP64 peak ratio are therefore attributable to the combined effect of higher arithmetic throughput, smaller memory traffic, improved cache behavior, and reduced nonlinear-solve work from the lower-precision implicit stages.  The effect is especially pronounced in comparisons involving the \(128/16\) pairing, since the full \(128/128\) baseline uses software quadruple precision, while the reduced-precision stage uses hardware-supported low precision arithmetic.  Finally, if future reduced-precision matrix kernels are routed through BF16 AMX rather than IEEE FP16 AVX-512, the compute roof can be much higher: Intel AMX uses two-dimensional tile registers and a tile matrix-multiplication unit, and Intel reports 1024 BF16 operations per cycle per core at full utilization \cite{IntelAMXBrief}, peek  $114.69~{\rm TFLOP/s}$. 
Thus the observed mixed-precision speedups reflect both hardware precision effects and algorithmic effects, rather than only the smaller size of the half-precision data type.

\section{Conclusions}\label{sec:conclusions}
In this work, we conducted a performance evaluation of mixed precision Runge--Kutta methods, comparing the original methods from~\cite{grant2022} with the newly proposed stabilized corrections in~\cite{driscoll2026}. This study extends the numerical experiments introduced in~\cite{burnett2021}, which focused on a simpler test problem. Our results confirm that runtime savings can also be obtained for more challenging test problems, such as Burgers' equation and the porous medium equation. Furthermore, this evaluation demonstrates that the new stabilized corrections not only improve the stability of the explicit corrections but also achieve significant runtime savings.

\section*{Acknowledgment}
This material is based upon work supported by 
AFOSR Grant No.\ FA9550-23-1-0037, 
Grant No.\ FA9550-24-1-0254, 
DOE Grant No.\ DE-SC0023164 and
DOE/NNSA Grant No. DE-NA0004265.
This project was initiated  while several authors were in residence at the
Institute for Computational and Experimental Research in Mathematics (ICERM)
during the program ``Empowering a Diverse Computational Mathematics Research
Community'' for which the authors acknowledge the National Science Foundation
under Grant No.\ DMS-1929284.

All computations for this project were performed on the Gautschi computing cluster at Purdue University.

Claude Opus 5 (Anthropic) was used to identify grammatical and punctuation errors in the manuscript and to assist in developing code subroutines for the numerical experiments. All text and code were reviewed and verified by the authors.

\section*{Author Contribution Statement}
CH provided the expertise on the mixed precision implementation,
wrote the mixed precision code, ran all the simulations, 
presented the data, and wrote the  first draft of the manuscript. 
ZJG provided the mixed precision framework.
ZJG and SG conceptualized the project, devised the correction approach, designed the experiments, 
and advised on the mixed precision implementation 
and understanding the stabilized corrections approach.
JD and TSK provided their expertise on designing the stabilization matrices $\Phi$ and 
on the problems shown, and assisted with coding, debugging, and understanding of the data.
AC provided the explanation of the computational speedups, allowing  a deep understanding of
the interaction between  the  hardware precision  and algorithmic effects on the chip used.

\bibliographystyle{IEEEtran}
\bibliography{references}

\end{document}